\documentclass{elsart}

\usepackage{enumerate}
\usepackage{amssymb}
\usepackage{url}

\usepackage{yjsco}
\usepackage{natbib}

\newcommand{\hilSym}{\mathop{\textrm{HF}}}

\newcommand{\La}{\mathcal{L}}
\newcommand{\s}{\ensuremath{^*}}
\newcommand{\p}{\ensuremath{^\prime}}

\newcommand{\setBuilder}[2]{\left\{#1|#2\right\}}
\newcommand{\hd}[1]{\textrm{in}_\leq(#1)}

\newcommand{\proofPart}[1]{{\bf $\bf{#1}$:}}

\newcommand{\Mi  }{(M1)}
\newcommand{\Mii }{(M2)}
\newcommand{\Miii}{(M3)}

\newcommand{\deSym}{\phi}
\newcommand{\de}[1]{\deSym\left(#1\right)}
\newcommand{\hil}[1]{\hilSym_{\hd{I_\La}}\left(#1\right)}
\newcommand{\ideal}[1]{\left\langle#1\right\rangle}

\begin{document}
\begin{frontmatter}

\author{Bjarke Hammersholt Roune}
\address{Department of Computer Science\\
University of Aarhus\\
IT-parken, Aabogade 34\\
DK-8200 Aarhus N, Denmark}
\ead[url]{\url{http://www.broune.com/}}

\title{Solving Thousand Digit Frobenius Problems Using Gr\"obner Bases}

\begin{abstract}
A Gr\"obner basis-based algorithm for solving the Frobenius Instance
Problem is presented, and this leads to an algorithm for solving the
Frobenius Problem that can handle numbers with thousands of
digits. Connections to irreducible decompositions and Hilbert
functions are also presented.
\end{abstract}
\end{frontmatter}

\section{Introduction}

Let $p_1,\ldots,p_n$ be relatively prime positive integers and let
$p=(p_1,\ldots,p_n)$. An integer is \emph{$p$-representable} if it can
be written as $v\cdot p$ for some $v\in\Nset^n$. Determining
$p$-representability is known as the \emph{Frobenius Instance
Problem}, and we present a new Gr\"obner basis-based algorithm that
solves it.\footnote{\cite{malkin} has discovered this independently.}
Our algorithm differs from the classical way to solve integer programs
using Gr\"obner bases due to \cite{IPGrobner} (see also
\cite{pottierToricGrobner}) by not adding any auxiliary variables to
the problem.

The \emph{Frobenius Number} $f_p\s$ is the largest integer that is not
$p$-representable, and such an integer exists by Proposition
\ref{frobExists} below. E.g. if $p=(6,10,15)$ then $f\s_p=29$.

The \emph{Frobenius Problem} is to compute the Frobenius Number
$f\s_p$. A recent algorithm due to \cite{frobPoint} can solve
Frobenius problems even if the $p_i$ have thousands of decimal
digits. Here we describe a simpler variant of that algorithm that
performs better.

The algorithm first computes a Gr\"obner basis and then determines the
Frobenius number based on that. \cite{frobPoint} use a novel algorithm
based on the Fundamental Domain that in effect computes this Gr\"obner
basis, and they report that this is much faster than using
Buchberger's algorithm in their implementation.

However, the benchmarks in section \ref{sec:benchmark} show that the
program 4ti2 [\cite{4ti2}], which implements Buchberger's algorithm in
a special case, outperforms the Fundamental Domain-based
implementation from \cite{frobPoint}. The record for random $p_i$ with
11 decimal digits was $n=11$, but we can now reach $n=13$. Performance
is also improved on smaller examples.

We show that the second step of the algorithm of \cite{frobPoint} can
be rephrased as the computation of the irreducible decomposition of a
certain monomial ideal. The final step of the algorithm is to maximize
a linear function over the decomposition, and we show that this
essentially computes the index of regularity of the ideal. We define
these terms when they are needed.

See \cite{mlfbFrobTest} for a geometric version of some of the ideas
in this paper in terms of maximal lattice free bodies. That paper is
joint work with Niels Lauritzen and Anders Nedergaard Jensen.

We refer to the book by \cite{frobBook} for more background on the
Frobenius Problem. We wish to thank Anders Nedergaard Jensen, Niels
Lauritzen, Daniel Lichtblau and Stan Wagon for helpful discussions
about Frobenius numbers.

\section{Preliminaries}

Let $e_i\in\Nset^n$, $i=1,\ldots,n$, be the vector whose entries are
all zero except that there is a 1 at position $i$. The
\emph{$p$-degree} of a vector $v\in\Zset^n$ is $v\cdot p$. If
$v\in\Zset^n$ then define $v^+\in\Nset^n$ as below and define
$v^-:=(-v)^+$.
\[
v^+_i:=\left\{
\begin{array}{cl}
v_i,&\textrm{for }v_i\geq0\\
0,&\textrm{for }v_i<0
\end{array}
\right.\]

We will need to consider the lattice ideal $I_\La$ defined by
\[ I_\La :=
\left\langle\left.x^{v^+}-x^{v^-}\right|v\in\Zset^n\textrm{ and }v\cdot p=0\right\rangle. \]

We will compute a Gr\"obner basis $G$ of $I_\La$. The term order
$\leq$ we will use first considers the $p$-degree of the exponent
vector of a term and then the reverse lexicographic order where
$x_1<x_2<\cdots<x_n$. 

\begin{exmp}
If $p=(2,3)$ then $1=x^{(0,0)}\leq x^{(1,0)} \leq x^{(3,0)}\leq
x^{(0,2)}\leq x^{(4,5)}=x_1^4x_2^5$.
\end{exmp}

We refer to \cite{cox} for more details.

\begin{prop}
\label{frobExists}
Only finitely many integers $t\in\Nset$ are not $p$-representable.
\end{prop}
\begin{pf}
As $p_1,\ldots,p_n$ are relatively prime, iterated use of the Extended
Euclidean Algorithm provides us with a vector $v\in\Zset^n$ such that
$v\cdot p=1$. Let $m:=\min_{i=1}^n v_i$ and define
$u:=v+p_1|m|(1,\ldots,1)$. Then $u+iv\in\Nset^n$ for $i=0,\ldots,p_1-1$
and therefore $(u+iv)\cdot p=u\cdot p+i$ are $p_1$ consecutive
$p$-representable numbers.
\end{pf}

\section{An Algorithm That Solves The Frobenius Instance Problem}
\label{sec:repTest}

We claim that the following algorithm determines if $t\in\Nset$ is
$p$-representable.

\begin{enumerate}[Step 1:]
\item Compute an $a\in\Zset^n$ such that $a\cdot p=t$, $a_1\leq 0$ and
$a_i\geq 0$ for $i=2,\ldots,n$.
\label{step1}

\item Divide $x^{a^+}-x^{a^-}$ by $G$ giving remainder
$x^w(x^{c^+}-x^{c^-})$ for some $w\in\Nset^n$ where $x^{c^-}\leq
x^{c^+}$.
\label{step2}

\item Then $t$ is $p$-representable if and only if $c\in\Nset^n$.
\label{step3}
\end{enumerate}

\begin{description}
\item[Step \ref{step1}]
We first find an $a\in\Zset^n$ such that $a\cdot p=t$. One way to do this
is to use the Extended Euclidean Algorithm iteratively to find a
$b\in\Zset^n$ such that $b\cdot p=1$ and then let $a:=tb$. As
$(-\sum_{i=2}^np_i,p_1,\ldots,p_1)$ has $p$-degree zero and the sign
pattern $(-,+,\cdots,+)$, we can assume that $a$ also has this sign
pattern by adding a sufficiently large multiple of this vector to $a$.

\item[Step \ref{step2}]
This step requires knowing the Gr\"obner basis $G$, and computing $G$
is the most time consuming part of the algorithm. Once $G$ has been
computed the polynomial division itself is comparatively fast.

\item[Step \ref{step3}]
Observe that the division algorithm ensures that $c$ will have no more
negative entries than $a$ does, so $c$ is negative at most in the
first coordinate. As $x^{c^+}$ is not reducible by $G$, the following
lemma tells us that if $t$ is $p$-representable then $c_1\geq 0$ so
that $c$ is a $p$-representation of $t$ since $c\cdot p=a\cdot p=t$.
\end{description}

\begin{lem}
\label{happyCamper}
Let $a\in\Zset^n$ such that $a_i\geq 0$ for $i=2,\ldots,n$ and
$a_1<0$. If $a\cdot p$ is $p$-representable then there exists a $g\in
G$ such that $\hd g|x^{a^+}$.
\end{lem}
\begin{pf}
As $a\cdot p$ is $p$-representable there exists a $b\in\Nset^n$ such that
$a\cdot p=b\cdot p$.  Letting $d:=a-b$ this implies that $d\cdot p=0$
whereby $h:=x^{d^+}-x^{d^-}\in I_\La$.  Thus there exists a $g\in G$
such that $\hd g|\hd h$. We will prove that $\hd g|\hd
h=x^{d^+}|x^{a^+}$.

\proofPart{\hd h=x^{d^+}} $d^+$ and $d^-$ have the same $p$-degree and
$d_1=a_1-b_1<0$.

\proofPart{x^{d^+}|x^{a^+}} This follows from $b\in\Nset^n$.
\end{pf} 

Note that $\leq$ can be replaced with any term order that first
considers the $p$-degree of the exponent vector and then the reverse
lexicographic order on the first variable.

\section{An Algorithm That Solves The Frobenius Problem}
\label{sec:frobFind}

The general idea of the algorithm is that we can represent $f\s_p$ by
a certain vector that has $p$-degree $f\s_p$, and that this vector has
certain properties (see Proposition \ref{algebraicMlfb}). It turns out
that only finitely many vectors have these properties, so we can look
through all of them, and then the one with maximal $p$-degree will be
the vector that represents $f\s_p$, i.e. it will have $p$-degree equal
to $f\s_p$ (see Proposition \ref{itWorks}).

Proposition \ref{algebraicMlfb} spells out the properties mentioned above.

\begin{prop}
\label{algebraicMlfb}
Let $c\in\Zset^n$ be the vector resulting from running the algorithm from
Section \ref{sec:repTest} on $t=f\s_p$. Then the following holds.
\begin{enumerate}
\item[\Mi] $c_1=-1$ and $c_i\geq 0$ for $i=2,\ldots,n$.
\item[\Mii] $x^{c^+}$ cannot be reduced by any element of $G$.
\item[\Miii] $x^{(c+e_i)^+}$ \emph{can} be reduced by some $g_i\in G$
for $i=2,\ldots,n$.
\end{enumerate}
\end{prop}
\begin{pf}
\proofPart{\textrm{\Mi}} It holds by construction that $c_1<0$ and that
$c_i\geq 0$ for $i=2,\ldots,n$. Let $c\p:=c+e_1$. Then
$x^{(c\p)^+}=x^{c^+}$ is not reducible by $G$ and $c\p\cdot p$ is
$p$-representable as $c\p\cdot p>f_p\s$. Then Lemma \ref{happyCamper}
implies that $c\p_1\geq 0$ whereby $c_1=-1$.

\proofPart{\textrm{\Mii}} 
This holds by construction.

\proofPart{\textrm{\Miii}} We see that $(c+e_i)\cdot p$ is $p$-representable as
 it is strictly larger than $f\s_p$. Thus Lemma \ref{happyCamper}
 provides a $g_i\in G$ such that $\hd{g_i}|x^{(c+e_i)^+}$.
\end{pf}

Let $M_p$ be the set of vectors $c\in\Zset^n$ that have the properties
\Mi{}, \Mii{} and \Miii{} from Proposition \ref{algebraicMlfb}. The
idea is to compute $M_p$ and then use Proposition \ref{itWorks} below
to find $f\s_p$.

\begin{prop}
\label{itWorks}
The following holds.
\begin{enumerate}[(i)]
\item $M_p$ is finite.
\item If $a\in M_p$ then $a\cdot p$ is not $p$-representable.
\item $f\s_p=\max\setBuilder{a\cdot p}{a\in M_p}$
\end{enumerate}
\end{prop}
\begin{pf}
\proofPart{\textrm{(i)}} Let $a\in M_p$ and $i\in\{2,\ldots,n\}$. Let
$g\in G$ such that $\hd g|x^{(a+e_i)^+}$. As $\hd g$ does not divide
$x^{a^+}$ we can infer that $a_i+1$ is the exponent of $x_i$ in $\hd
g$. Thus there are at most $|G|$ possibilities for what $a_i$ can be.

\proofPart{\textrm{(ii)}} Lemma \ref{happyCamper} shows that this
follows from \Mi{} and \Mii{}.

\proofPart{\textrm{(iii)}} Proposition \ref{algebraicMlfb} shows that there is an
$a\in M_p$ such that $a\cdot p=f\s_p$. This and part (ii) above shows
what we need.
\end{pf}

We can compute $M_p$ as follows. We saw in the proof of Proposition
\ref{itWorks} that if $i\in\{2,\ldots,n\}$ and $a\in M_p$ then $a_i+1$
is the exponent of $x_i$ in $\hd{g_i}$ for some $g_i\in G$. Thus we
can run through all possible values of $a_i$ for $i=2,\ldots,n$ and
only keep those $a$ that have properties \Mi{}, \Mii{} and \Miii{}.

This algorithm is easy to understand and implement, but it requires us
to look through up to ${|G|}^{n-1}$ possibilities. The External Corner
Algorithm from \cite{frobPoint} is a more efficient algorithm for
computing $M_p$ which usually dramatically reduces the number of
possibilities that need to be examined.

In Section \ref{sec:decomConnection} we will define the irreducible
decomposition of a monomial ideal, and we will prove that $M_p$
corresponds to the irreducible decomposition of the initial ideal of
$I_{\La}$. \cite{stairIrr} shows how an algorithm that is very similar
to the External Corner Algorithm can compute irreducible
decompositions of monomial ideals in general in much less time than
the best competing programs.

\section{A Connection To Monomial Irreducible Decompositions}
\label{sec:decomConnection}

We need a few more definitions. Let $I$ be a monomial ideal and define
the function $\deSym$ by $\de v=\left\langle
x_i^{v_i}|v_i>0\right\rangle$ for $v\in\Nset^n$ except that $\de
1=\ideal 1$. An ideal of the form $\de v$ is called
\emph{irreducible} and the \emph{irredundant irreducible
decomposition} of $I$ is the unique minimal subset $D\subseteq\Nset^n$
such that $I=\bigcap_{v\in D}\de v$. Thus the irredundant irreducible
decomposition of $I:=\ideal{x_1^2,x_1x_2}$ is $\{(1,0),(2,1)\}$ as
$I=\ideal{x_1}\cap\ideal{x_1^2,x_2}$.

An ideal is \emph{artinian} if there exists a $t\in\Nset$ such that
$x_i^t\in I$ for $i=1,\ldots,n$. Note that $\hd{I_\La}$ is artinian if
we first project out the variable $x_1$, since $x_i^{p_1}-x_1^{p_i}\in
I_\La$ and therefore $\hd{x_i^{p_1}-x_1^{p_i}}=x_i^{p_1}\in\hd{I_\La}$
for $i=2,\ldots,n$.

We claimed in Section \ref{sec:frobFind} that $M_p$ corresponds to the
irreducible decomposition of the initial ideal of $I_\La$, and
Proposition \ref{decom} proves this claim.

\begin{prop}
\label{decom}
The set $M\p:=\{(a_1+1,\ldots,a_n+1)|a\in M_p\}$ is the irreducible
irredundant decomposition $D$ of $\hd{I_\La}$.
\end{prop}
\begin{pf}
Let $a\in\Zset^n$ and let $a\p:=a-\sum_{i=1}^ne_i$. Consider the
following statements.
\begin{enumerate}
\item $a\in D$
\label{s1}
\item $a_1=\phantom{-}0$, $a_i\geq 1$, $\hd{I_\La}\subseteq\de{a}$ and
$\hd{I_\La}\not\subseteq\de{a+e_i}$ for $i=2,\ldots,n$.
\label{s2}
\item $a\p_1=-1$, $a\p_i\geq 0$, $x^{(a\p)^+}\notin \hd{I_\La}$ and
$x^{(a\p+e_i)^+}\in \hd{I_\La}$ for $i=2,\ldots,n$.
\label{s3}

\item $a\p\in M_p$
\label{s4}
\end{enumerate}
We will prove that these statements are equivalent.

\proofPart{\textrm{(\ref{s1})}\Leftrightarrow\textrm{(\ref{s2})}} The
initial ideal of $I_\La$ contains no monomial that is divisible by
$x_1$, and this implies that $a_1=0$ for $a\in D$. Now that we have
handled the first entry, we can project out $x_1$ and thereby get an
artinian ideal.

When working with artinian ideals the elements of the decomposition
consists of vectors without zero entries, and if $a$ has no zero
entries, then $\de {a+e_i}\subsetneq\de a$.

To get the irredundant irreducible decomposition of an ideal we write
the ideal as an intersection of irreducible ideals that are as small
as possible, and this is exactly what (2) expresses since the
projected ideal is artinian.

\proofPart{\textrm{(\ref{s2})}\Leftrightarrow\textrm{(\ref{s3})}}
By Lemma \ref{tech} below.

\proofPart{\textrm{(\ref{s3})}\Leftrightarrow\textrm{(\ref{s4})}}
By definition.
\end{pf}

\begin{lem}
\label{tech}
Let $m\in\Nset^n$ and let $I$ be a non-zero monomial ideal. Then $x^m\in
I$ if and only if $I\not\subseteq\de{m+\sum_{i=1}^n e_i}$.
\end{lem}
\begin{pf}
Let $a\in I$ be a monomial. Then $a|x^m$ if and only if
$a\notin\de{m+\sum_{i=1}^ne_i}$.
\end{pf}

\section{A Connection To The Hilbert Function}

The \emph{$p$-weighted Hilbert function} $\hilSym_I\colon\Nset\to\Nset$ of a
monomial ideal $I$ is defined such that $\hilSym_I(t)$ is the number
of monomials $x^m\notin I$ where $m$ has $p$-degree $t$.

\begin{prop}
$\hil t$ is equal to 1 if $t$ is $p$-representable and 0 otherwise.
\end{prop}
\begin{pf}
\proofPart{\hil t\leq 1} Let $m_1,m_2$ be two vectors of $p$-degree
$t$. Then $x^{m_1}-x^{m_2}\in I_{\La}$ whereby the initial term is in
$\hd{I_\La}$. Thus at most one of $x^{m_1}$ and $x^{m_2}$ is not in
$\hd{I_\La}$.

\proofPart{\hil t=1\Rightarrow \textrm{representability}} If $\hil t=1$ then
there is a monomial $x^m\notin\hd{I_{\La}}$ where $m$ has $p$-degree
$t$. Thus $t$ is $p$-representable.

\proofPart{\textrm{representability}\Rightarrow\hil t=1} Run the
algorithm from Section \ref{sec:repTest} on $t$. The resulting vector
$c\in\Nset^n$ is such that $x^c\notin\hd{I_\La$}.
\end{pf}

We can infer that $\hil{f\s_p}=0$ and that $\hil{t}=1$ for all
integers $t>f\s_p$. The integer at which the Hilbert function becomes
equal to a polynomial is known as the \emph{index of regularity}, so
in this case the index of regularity is $f\s_p+1$.

The elements of $M_p$ correspond to the maximal monomials outside of
$\hd{I_\La}$ according to divisibility (disregarding the first
variable), and what the algorithm from Section \ref{sec:frobFind} does
is to maximize the dot product with $p$ over $M_p$. This amounts to
maximizing the dot product over the vectors $m$ such that
$x^{m+e_1}\notin \hd{I_\La}$ and $m_1=-1$.

We can multiply any monomial not in $\hd{I_\La}$ with $x_1$ and get
something still not in $\hd{I_\La}$, so anything that goes on in the
first variable does not prevent the Hilbert function from becoming a
polynomial. Thus the algorithm can be interpreted as finding the
maximal point that prevents the Hilbert function from becoming a
polynomial, which is to say that it computes the index of regularity.

We conclude that the algorithm of \cite{frobPoint} can be interpreted
as computing an index of regularity.

\section{Benchmarks}
\label{sec:benchmark}

\cite{frobby} has written an implementation called Frobby of the
algorithm described in this paper, and here we compare Frobby to the
implementation of \cite{frobPoint}. Figure \ref{data} displays the
collected data. ``Intractable'' means intractable according to the
authors of that software package. Note that the most time consuming
part of the algorithm usually is to compute the Gr\"obner
basis. Frobby and Mathematica are the only programs that can handle
input numbers $p_i$ as large as those in figure \ref{data}.

Frobby uses the program fplll due to \cite{fplll} to obtain an
LLL-reduced lattice basis and then computes the Gr\"obner basis of the
lattice ideal $I_\La$ from that using the program 4ti2
[\cite{4ti2}]. Frobby then computes the Frobenius number by computing
an irreducible decomposition of $\hd{I_\La}$ using the algorithm due
to \cite{stairIrr}, which is similar to the External Corner
Algorithm. Frobby is written in C++ and the source code is available
under the GNU General Public License (GPL).

The implementation of \cite{frobPoint} is available as a part of
Mathematica and is written mostly in C. The major difference from
Frobby is that the Gr\"obner basis is computed using a different
algorithm than Buchberger's.

All the inputs were randomly generated using genuinely random
radioactive decay via the service provided by \cite{hotbits} except
the $n=11$ input which was provided by Stan Wagon who pseudo-randomly
generated it using Mathematica. We wish to thank Daniel Lichtblau for
carrying out the Mathematica benchmarks.

All the benchmarks were run on machines with a 3.0 GHz Pentium 4 CPU
with 1 GB RAM except the Mathematica benchmark on the $n=11$ input
which was run on a 3.2 GHz Pentium 4.

\begin{figure}
\centering
\begin{tabular}{|c|c|c|c|c|c|}
\hline
$n$&$\lfloor\log_{10}(\min p_i)\rfloor+1$&
 $|G|$& Mathematica&Frobby\\
\hline\hline
4  & 10000 & 7    & 620s  & 3.5s \\

4  & 800 & 10     & 4.0s  & 0.3s \\
5  & 150 & 34     & 3.8s  & 0.3s  \\
6  & 70  & 131    & 5.5s  & 0.3s \\
6  & 80  & 148    & 5.6s  & 0.3s \\
6  & 90  & 112    & 6.8s  & 0.3s \\
6  & 100 & 140    & 9.5s  & 0.3s \\
8  & 30  & 2099   & 80.2s & 11.8s \\

11 & 11  & 27037  & 43.3h & 0.5h \\

12 & 11  & 56693  & intractable & 2.5h \\
13 & 11  & 170835 & intractable & 49.7h  \\
\hline
\end{tabular}
\caption{The benchmark data.}
\label{data}
\end{figure}

\bibliography{references}
\bibliographystyle{elsart-harv}

\end{document}